\documentclass[12pt]{article}

\usepackage[T2A]{fontenc}
\usepackage[utf8]{inputenc}

\usepackage{amsfonts}
\usepackage{amssymb}
\usepackage{amsmath}
\usepackage{amsthm}

\def \le {\leqslant}
\def \ge {\geqslant}

\begin{document}

\begin{Large}
\centerline{\bf Eine Bemerkung  \"uber  positiv definite  quadratische Formen }
\centerline{\bf und  rationale Punkte}
 
\end{Large}
\vskip+1cm
\begin{large}
\centerline{von Nikolay Moshchevitin\footnote{ unterst\"{u}tzt durch
RFBR 15-01-05700а.
 } (Moskau)}
 \end{large}
\vskip+2cm

{\bf 1. Resultate \"uber  rationale Punkte.}
  
  Sei 
  $$
 f({\bf x}) = \sum_{i,j =1}^{n} f_{i,j} x_i x_j,\,\,\,\,\,\,
 f_{i,j}  = f_{j,i},\,\,\,\,\,\, 
   {\bf x} = (x_1,...,x_n) \in \mathbb{R}^n
 $$
 eine positiv definite quadratische Form   mit ganzen Koeffizienten.
    F\"ur die Form $f$  betrachten wir eine
   indefinite Form
  $$
 F({\bf z}) =  f({\bf x}) - y^2
 $$
in  $n+1$ Variablen
 $ {\bf z} = (x_1,...,x_n, y) \in \mathbb{R}^{n+1}$
 und  den K\"orper
  $$
  \hbox{\got S}_f =
   \{ {\bf x}\in \mathbb{R}^n: |f( {\bf x} )| < 1\},
   $$ 
   dessen Volumen wir durch $\hbox{\got v}_f$ bezeichnen.
   
   Wir erinnern uns, dass eine quadratische Form $Q ({\bf z})$    isotrop hei\ss t, wenn
 ein Gitterpunkt ${\bf g} \neq{\bf0}$ mit $Q ({\bf g})=0$ existiert. 
    
    In der vorliegenden Note
    werden wir     
    folgende Behauptung zeigen.

    {\bf Satz 1.}\,\,{\it
     Sei 
     $$f(\pmb\alpha ) =
     f(\alpha_1,...,\alpha_n) =1.
     $$
     Definieren wir die Konstante
    $$
    \kappa_f = 
    6\max\left\{ 1,
     \frac{(n+1)^{2(n+1)} \cdot 6^{2n}\cdot 2^{2n^2}}{\hbox{\got v}_f^{2(n+1)}}
    \right\}
    $$
   Falls die Form  $F({\bf z})$
   isotrop ist, gibt es
    f\"ur jedes $ T\ge 1$  einen  ratiolalen Punkt
    $$
    {\bf r} = \left(\frac{a_1}{q},...,\frac{a_n}{q}\right)
    ,
    \,\,\,
    a_1,...,a_n \in \mathbb{Z},\,\,\, q\in \mathbb{Z},
    $$
    mit den folgenden Eigenschaften:  
   $$
   1\le q\le T,
   $$
   \begin{equation}\label{qqa}
   f\left(\pmb\alpha -{\bf r}\right)=
   f\left(\alpha_1-\frac{a_1}{q},...,\alpha_n-\frac{a_n}{q}\right)
   \le   \frac{\kappa_f}{qT}
   \end{equation}
   und
   $$
   f({\bf r}) =
      f\left(\frac{a_1}{q},...,\frac{a_n}{q}\right) =1.
   $$
   
   }
   
 Diese Behauptung ist eine effektive  Version eines h\"ubschen Satzes von Fishman, Kleinbock, Merrill und Simmons [7]. 
  Kleinbock und Merrill [6] hatten vor  kurzem das Resultat des Satzes 1
  f\"ur die einfachste positive Form
 \begin{equation}\label{sero}
 f_0 ({\bf x}) =  x_1^2+...+x_n^2
 \end{equation}
 bewiesen,
ohne eine explizite Konstante in (\ref{qqa}).
Ein allgemeines Ergebnis war später  von
Fishman,
 Kleinbock, Merrill und Simmons [6] bewiesen werden.
Sie hatten keine explizite Konstante auch und sie benutzten die Methoden der Theorie 
 Dynamischer Systeme.

  Hier  geben wir einen kurzen Beweis des Satzes 1, der  nur
  den
   Satz von Minkowski \"uber die sukzessiven Minima und einen Nullstellensatz 
 f\"ur die isotropen Form  von Cassels [3] verwendet.  
 Dieser Nullstellensatz  wurde 1955 von Cassels  [1] bewiesen (sehe auch die Arbeit von Davenport [4]).
 Birch und Davenport [1], Schlickewei [9], 
 Schulze-Pillot [12], 
 Schlickewei und Schmidt [10,11]
 haben dieses Resultat verallgemeinert.
 Eine interessante \"Ubersicht wurde von Fukshansky [5] geschrieben.
  Ein kurzer Beweis  des Satzes 1 f\"ur $f_0$ in dem Fall $ n=3$
 wurde in [8] gegeben.
 
 \"Altere Resultate \"uber
 die Aproximation von Punkten auf der  quadratischen Fl\"ache
 durch rationale Punkte  der Fl\"ache
  sind in dem Buch von Cassels ([3], Kap. 6, \S \,8, 9) besprochen.
  
  \vskip+0.5cm
 
 {\bf 2. Quadratische Forme und Automorphismen.}

 Wir brauchen
 die einfachste  Form
(\ref{sero})  
 und die entsprechende  Form
 $$
 F_0({\bf z} ) =  x_1^2+...+x_n^2 - y^2.
 $$
  Wir betrachten die Kugel
 $$
    \hbox{\got S}
  =\{ {\bf x}:\,\,
  x_1^2+...+x_n^2< 1 \}\subset \mathbb{R}^n
  $$ 
  mit dem Volumen
  $
  \hbox{\got o}_n = {\rm vol} \hbox{\got S} .$

 Wir brauchen auch die indefinite Form
  $$
 F_1({\bf w} ) =  x_1^2+...+x_{n-1}^2 + \xi\eta
 $$
 in  $n+1$ Variablen
 $ {\bf w} = (x_1,...,x_{n-1},\xi, \eta) \in \mathbb{R}^{n+1}$.
 Es ist klar, dass
 $$
   F_0({\bf z} ) =F_1({\cal B} {\bf z} )
   $$
 mit
  $$
  {\cal B}=
  \left(\begin{array}{cc cccc}
  1&\cdots&0 &0&0\cr
  \vdots&\vdots&\vdots&\vdots&\vdots\cr
    0&\cdots &1&0& 0\cr
  0&\cdots &0&1&- 1\cr 
    0&\cdots&0&1&1
     \end{array}
\right).
  $$
  Die Form $F_0$ hat  eine Familie der Automorphismen
  $$
{\cal D}_t =
\left(
\begin{array}{ccccccc}
1 &0&\cdots &0&0&0\cr
0&1&\cdots &0&0&0\cr
\vdots&\vdots&\vdots&\vdots&\vdots&\vdots\cr
0&0&\cdots& 1&0&0\cr
0&0&\cdots&0&t^{-1}&0\cr
0&0&\cdots &0&0&t
\end{array}
\right),
\,\,\,\,
    t\in \mathbb{R}_+, 
$$
so dass 
 $$
{\cal G}_t =
{\cal B}^{-1}
{\cal D}_t{\cal B} =
\left(
\begin{array}{ccccccc}
1 &0&\cdots &0&0&0\cr
0&1&\cdots &0&0&0\cr
\vdots&\vdots&\vdots&\vdots&\vdots&\vdots\cr
0&0&\cdots& 1&0&0\cr
0&0&\cdots&0&
\frac{1}{2}\left(t+\frac{1}{t}\right)
 &
 \frac{1}{2}\left(t-\frac{1}{t}\right)
 \cr
0&0&\cdots &0&\frac{1}{2}\left(t-\frac{1}{t}\right)&
\frac{1}{2}\left(t+\frac{1}{t}\right)
\end{array}
\right),
\,\,\,\,,
    $$
   Automorphismen der Form $F_0$ sind:
    $$
    F_0 ({\cal G}_t {\bf z}) =
   F_0 ( {\bf z})
 ,\,\,\,\,\,
  \forall\,t\,\,\forall\, {\bf z}
 .
 $$
 F\"ur  einen  Punkt $\pmb\beta \in \hbox{\got S}$ k\"onnen wir eine orthogonale Matrix
 $R_{\pmb\beta}$
 finden, f\"ur die gilt
 $$
 \pmb\beta = R_{\pmb\beta} {\bf e},
 $$
 wo ${\bf e} = (\underbrace{0,...,0}_{n-1},1)^T$. Es ist klar, dass $R_{\pmb\beta}$  ein ``abstandserhaltenden Automorphismus'' der Form $F_0$ ist, der
 den euklidischen Abstand zwischen Punkten bewahrt.
 
 Wir definieren die Punkte
 $$
 \overline{\bf e} = 
 \left(
 \begin{array}{c}
  0\cr
  \vdots\cr
  0\cr
  1
  \cr
  1
 \end{array}
\right),
 \,\,
 \overline{\pmb\alpha} = 
 \left(
 \begin{array}{c}
  \alpha_1\cr
  \vdots\cr
  \alpha_{n-1}\cr
  \alpha_n\cr
 1
 \end{array}
\right),
 \,\,
 \overline{\pmb\beta} = 
 \left(
 \begin{array}{c}
  \beta_1\cr
  \vdots\cr
  \beta_{n-1`}\cr
  \beta_n\cr
  1
 \end{array}
\right)
\in\mathbb{R}^{n+1}
$$
und die Matrizen
  $$
 {\cal W}_f
  =
   \left(
  \begin{array}{cc}
   W_{f\,\, n\times n} &  0_{1\times n}\cr
   0_{n\times 1} & 1 
  \end{array}
\right)
,\,\,\,\,\,
 {\cal R}_{\pmb\beta}
  =
   \left(
  \begin{array}{cc}
   R_{\pmb\beta\,\, n\times n} &  0_{1\times n}\cr
   0_{n\times 1} & 1 
  \end{array}
\right).
    $$
 Dann
 \begin{equation}\label{q1}
  \overline{\pmb\alpha} = {\cal W}_f\overline{\pmb\beta},
  \end{equation}
  \begin{equation}\label{q2}
  \overline{\pmb\beta}  ={\cal R}_{\pmb\beta} \overline{\bf e}
  \end{equation}
  und
  \begin{equation}\label{q3}
 t\cdot \overline{\pmb\beta}  ={\cal R}_{\pmb\beta} {\cal G}_t \overline{\bf e} 
 .
 \end{equation}

 Wir betrachten den Zylinder
   $$
 \hbox{\got K} =
\left\{ {\bf z} = (x_1.,,,x_n,y) \in \mathbb{R}^{n+1}:\,\,
|y|<1,\,
x_1^2+...+x_{n-1}^2+x_n^2< 1 \right\}
$$
mit dem Volumen $ {\rm vol}\, \hbox{\got K} = 2\hbox{\got o}_n$.
Sei $   \hbox{\got K}_f(\pmb\alpha, t)  =
{\cal W} {\cal R}_{\pmb\beta (\pmb\alpha)} {\cal G}_t \hbox{\got K}$.
Wegen  
$$
 \hbox{\got K}\subset \{ {\bf z}: \,\, |F_0({\bf z}) |< 1\},
 $$
 \begin{equation}\label{ffff}
 F_0 ({\cal R}_\beta {\bf z}) =  F_0 ({\cal G}_t {\bf z}) =F_0 ( {\bf z}),\,\,\,
 F(W{\bf x}) = F_0({\bf x}),
 \end{equation}
 haben wir
 \begin{equation}\label{kaa}
 \hbox{\got K}_f (\pmb\alpha , t)\subset \{ {\bf z}: \,\, |F({\bf z}) |< 1\}.
  \end{equation}
  
  Wir bemerken, dass
 \begin{equation}\label{wwo}
  {\rm vol}\,  \hbox{\got K}_f (\pmb\alpha , t) =
  |{\rm det}\, W_f\,\cdot
 \,{\rm det}\, R_{\pmb\beta (\pmb\alpha)}\,\cdot\,
 {\rm det }\, {\cal G}_t\,\,\cdot 2\hbox{\got o}_n | = 
 |{\rm det}\, W_f |\,\cdot 2\hbox{\got o}_n = 2\hbox{\got v}_f.
  \end{equation}

   \vskip+0.5cm
 {\bf 3. Sukzessive Minima.}

 Wir betrachten die sukzessiven Minima  $\lambda_j, \, 1\le j \le n+1$ des K\"orpers
$\hbox{\got K}_f (\pmb\alpha , t)$.
Aus dem Satz von  Minkowski folgt
$$
\lambda_1^n\lambda_{n+1} \le \lambda_1\cdots \lambda_{n+1}\le \frac{ 2^{n+1}}{ 
{\rm vol} 
\,\hbox{\got K}_f (\pmb\alpha , t)} ,
$$
sodass
\begin{equation}\label{le}
\lambda_{n+1} \le 
\frac{ 2^{n+1}}{ \lambda_1^n\cdot
{\rm vol} 
\,\hbox{\got K}_f (\pmb\alpha , t)} =
\frac{ 2^{n}}{ \lambda_1^n\cdot
\hbox{\got v}_f}.
\end{equation}

Wir unterscheiden zwei Fälle.

{\bf  Fall 1.} Sei $\lambda_1 < 1$. \,\,\,
In diesem Fall existiert ein Gitterpunkt
$
{\bf g}\in \hbox{\got K}_f (\pmb\alpha , t) 
,\,\,
{\bf g}\neq{\bf 0}.
 $
 Aber $F({\bf g})$ ganzzahlig ist, sodass
 aus (\ref{kaa})  folgt $ F({\bf g})=0$.
 
 {\bf  Fall 2.} Sei $\lambda_1 \ge 1$. \,\,\, Dann aus 
  (\ref{le})
    folgt
  \begin{equation}\label{n1}
  \lambda_{n+1} \le \frac{ 2^{n}}{ \hbox{\got v}_f} .
  \end{equation}
  In diesem Fall existieren unabh\"angig Gitterpunkte
  \begin{equation}\label{lll}
  {\bf g}_1,...,{\bf g}_{n+1} \in   \lambda_{n+1}
  \overline{
  \hbox{\got K}_f (\pmb\alpha , t)
}.
\end{equation}
Die Form $F$ ist isotrop und darum die Form
$$
 Q({\pmb\xi }) =  F(\xi_1{\bf g}_1+...+\xi_{n+1}{\bf g}_{n+1} ) = 
\sum_{i,j =1}^{n+1} Q_{i,j} \xi_i\xi_j,
 \,\,\, \pmb\xi =(\xi_1,...,\xi_{n+1})\in \mathbb{Z}^{n+1}
 $$
  ist isotrop auch.
 Dann  Hilfssatz 
8.1 aus Kap. 6 [3] gibt einen Gitterpunkt $\pmb\zeta=(\zeta_1,...,\zeta_{n+1})\in \mathbb{Z}^{n+1}\setminus \{{\bf 0}\}$ mit
$
 Q(\pmb\zeta) =0
$
und
\begin{equation}\label{ll}
 \max_{1\le j \le n+1} |\zeta_j| \le \left(3\sum_{1\le i,j\le n+1} |Q_{i,j}|\right)^{n/2}.
\end{equation}
Aber
  $$
 Q_{j,j} = F({\bf g}_j),
 \,\,\,
 Q_{i,j}
 =\frac{1}{2}\left( F({\bf g}_i+{\bf g}_j)
 - F({\bf g}_j)- F({\bf g}_i)
 \right) , i \neq j
 $$
 Aus (\ref{lll})  folgt  $   {\bf g}_i+
 {\bf g}_j\in 2  \lambda_{n+1}
  \overline{
  \hbox{\got K}_f (\pmb\alpha , t)
}$, sodass liefern
 ({\ref{lll}) und (\ref{kaa}) 
 $$
|F( {\bf g}_j)| \le \lambda_{n+1}^2 ,\,\,\,
|F(  {\bf g}_i+
 {\bf g}_j)| \le 4\lambda_{n+1}^2, 
 $$
und darum
\begin{equation}\label{llll}
 \max_{1\le i,j\le n+1} |Q_{i,j}|\le 3\lambda_{n+1}^2.
\end{equation}

Sei
$${\bf g} =
\zeta_1{\bf g}_1+...+\zeta_{n+1}{\bf g}_{n+1}
  .$$
Dann $F({\bf g}) =0$ und
aus (\ref{n1}), (\ref{lll}), (\ref{ll}) und (\ref{llll}) folgt
$$
{\bf g} \in
(n+1)^{n+1} 3^n\lambda_{n+1}^{n+1}
  \hbox{\got K}_f (\pmb\alpha , t)
  \subset
  \frac{(n+1)^{n+1} \cdot 6^n\cdot 2^{n^2}}{\hbox{\got v}_f^{n+1}} \, \hbox{\got K}_f (\pmb\alpha , t) .  
  $$

  In beiden Fällen
  existiert  ein Gitterpunkt
  \begin{equation}\label{g}
  {\bf g} = (a_1,...,a_n,q )\in C_f
  \, \hbox{\got K}_f (\pmb\alpha , t) \setminus \{{\bf 0}\},
  \,\,\,\,
  C_f = \max \left\{ 1, 
  \frac{(n+1)^{n+1} \cdot 6^n\cdot 2^{n^2}}{\hbox{\got v}_f^{n+1}}\right\}   
  \end{equation}
  mit
  $F({\bf g}) =0$.

   \newpage

  {\bf 4. Beweis des Satzes 1.}
 
 Sei
 $$
 t =\frac{2T}{3C_f.
 }
 $$  
Betrachten wir den Gitterpunkt ${\bf g} $ aus (\ref{g}). Wir definieren die Punkte 
\begin{equation}\label{p1}
{\bf u} =(u_1,...,u_n, q) = {\cal W}_f^{-1} {\bf g} \in {\cal R}_{\pmb\beta (\pmb\alpha)} {\cal G}_t C_f\hbox{\got K},\,\,\,\,\,
\,\,\,\,\,\,\,\,\,\,\,\,\,\,\,\,\,\,\,\,\,\,\,\,\,\,\,\,\,\,\,\,\,\,\,\,\,\,\,\,\,\,\,\,\,\,\,\,\,\,\,\,\,\,\,\,\,\,
\,\,\,\,\,\,\,\,\,\,\,\,\,\,\,\,\,\,\,\,\,\,\,\,
{\bf g} ={\cal W}_f{\bf u};
\end{equation}
\begin{equation}\label{p2}
{\bf v} = (v_1,...,v_n, q) =  {\cal R}_{\pmb\beta (\pmb\alpha)}^{-1}  {\bf  u} =
{\cal R}_{\pmb\beta (\pmb\alpha)}^{-1}  {\cal W}_f^{-1} {\bf g}
  \in  {\cal G}_t C_f\hbox{\got K},   \,\,\,\,\,\,\,\,\,\,\,\,\,\,\,\,\,\,\,\,
  \,\,\,\,\,\,\,\,\,\,\,\,\,\,\,\,\,\,\,\,\,\,\,\,\,\,\,\,\,\,\,\,\,\,\,\,
    {\bf u } ={\cal R}_{\pmb\beta (\pmb\alpha )}{\bf v};
\end{equation}
\begin{equation}\label{p3}
{\bf w} = (w_1,...,w_n, w_{n+1}) = {\cal G}_t^{-1}{\bf v} = {\cal G}_t^{-1} {\cal R}_{\pmb\beta (\pmb\alpha)}^{-1}  {\bf  u} =
{\cal G}_t^{-1} 
{\cal R}_{\pmb\beta (\pmb\alpha)}^{-1}  {\cal W}_f^{-1} {\bf g}
  \in C_f \hbox{\got K},\,\,\,\,\,\,\,\,\,\,\,\,\,\,\, {\bf v} ={\cal G}_t{\bf w};
\end{equation}
so
$$
w_j = v_j ,\,\,\, 1\le j \le n-1,
$$
$$
w_n =\frac{1}{2}
\left( \frac{1}{t}+t\right)v_n + 
\frac{1}{2}\left( \frac{1}{t}-t\right)q,
$$
$$
w_{n+1} =\frac{1}{2}
\left( \frac{1}{t}-t\right)v_n + 
\frac{1}{2}\left( \frac{1}{t}+t\right)q
$$
und
$$
F_0 ({\bf w}) =
F_0 ({\bf v}) =
v_1^2+...+v_{n-1}^2+v_n^2 - q^2 =
$$
$$
=
v_0^2+...+v_{n-1}^2 
+
\left( \frac{1}{2}
\left( \frac{1}{t}+t\right)v_n + 
\frac{1}{2}\left( \frac{1}{t}-t\right)q\right)^2-
\left( \frac{1}{2}
\left( \frac{1}{t}-t\right)v_n + 
\frac{1}{2}\left( \frac{1}{t}+t\right)q
\right)^2.
$$
Wegen ${\bf w} \in C_f\hbox{\got K}$, gilt
$|w_n|,|w_{n+1}|<C_f$ und
$
|w_n+w_{n+1}|<2C_f,\, |w_n - w_{n+1}|<2C_f
$.
Daher ist
\begin{equation}\label{minu}
 |q-v_n|< \frac{2C_f}{t},
\end{equation}
und
\begin{equation}\label{pilu}
 |q+v_n|< {2C_f}{t},
\end{equation}
Aus (\ref{minu}), (\ref{pilu}) 
folgt
$$
2q = q+v_n+q-v_n < \frac{2C_f}{t}+2C_ft < 3C_ft
,$$
oder
\begin{equation}\label{qu}
 q\le \frac{3C_ft}{2} =T.
\end{equation}
 Aus (\ref{q1}), (\ref{p1}), (\ref{q2}), (\ref{p2}) und (\ref{ffff})  folgt
 $$
 F({\bf g} -q\overline{\pmb\alpha})
 =
 F({\cal W}_f ({\bf u} -q\overline{\pmb\beta}))=
 F_0({\bf u} -q\overline{\pmb\beta})
 = F_0( {\rm \cal R}_{\pmb\beta (\pmb\alpha)}({\bf v} -q\overline{\bf e}))=
  F_0({\bf v} -q\overline{\bf e})
=
v_1^2+...+v_{n-1}^2 +
(v_n-q)^2. $$
  Aber $ v_1^2+...+v_{n-1}^2 +v_n^2 - q^2 =F_0({\bf v}) = F({\bf g}) = 0$. Daraus folgt
  \begin{equation}\label{finale}
 |F(
 {\bf g} -q\overline{\pmb\alpha})
  |=
|v_1^2+...+v_{n-1}^2 +
(v_n-q)^2  
 |=
 |q^2 -v_n^2 + (v_n - q)^2| =
 |2q(q-v_n)| <\frac{4C_fq}{t}.
 \end{equation}
 (In der letzten Ungleichung verwenden wir
(\ref{minu}).)  
 F\"ur den rationalen Punkt ${\bf r} =\left(\frac{a_1}{q},...,\frac{a_n}{q}\right) \in \mathbb{Q}^n$
   folgt
  \begin{equation}\label{finale}
  f\left(\pmb\alpha -{\bf r}\right) \le   \frac{6C_f^2}{qT}.
   \end{equation}
   Nun haben wir (\ref{qu}) und (\ref{finale}), und Satz 1 bewiesen ist.

 {\bf 5.  Linear unabh\"angige  L\"osungen.}

 Das folgende Ergebnis wurde 1983 von Schulze-Pillot [12] gezeigt
  (seht auch Folgerung 1 aus [11]):

 {\it
 Sei 
  $$
  Q({\pmb\xi }) =
  \sum_{i,j =1}^{m} Q_{i,j} \xi_i\xi_j
    , \,\,\, \pmb\xi =(\xi_1,...,\xi_{m})\in \mathbb{R}^{m}
 $$
 eine isotrope  quadratische Form mit ganzrationalen Koeffizienten.
 Dann gibt es linear unabh\"angige
ganzzahlige Punkte 
${\pmb\zeta}^k = (\zeta_1^k,...,\zeta_m^k) ,\,\,
1\le k \le m$
mit 
$
Q({\pmb\zeta }^k) =0$ und mit
$$
\prod_{k=1}^m \max_i |\zeta_i^k| \le c_m \left( \max_{i.j} |Q_{i.j}|
\right)^{\frac{m^2}{2}-1}
,$$
wobei $c_m >0$ eine effektive Konstante ist. 
 }

 Aus diesem Ergebnis folgt:

 {\bf Satz 2.}\,\,{\it
     Sei 
     $
     f(\alpha_1,...,\alpha_n) =1.
     $
       Falls die Form  $F({\bf z})$
   isotrop ist, gibt es
eine effektive Konstante $\kappa_f^*>0$
   mit den folgenden Eigenschaften.
     F\"ur jedes $ T\ge 1$   gibt es $(n+1)$
 unabh\"angige
ganzzahlige Punkte 
$$
(a_1^k,...,a_n^k, q^k) \in \mathbb{Z}^{n+1},\,\,\,\, 1\le k \le n+1
$$
mit
$   1\le q^k\le T$,
$   f\left(\alpha_1-\frac{a_1^k}{q^k},...,\alpha_n-\frac{a_n^k}{q^k}\right)
   \le   \frac{\kappa_f^*}{q^kT}
$
   und
   $
        f\left(\frac{a_1^k}{q^k},...,\frac{a_n^k}{q^k}\right) =1.
   $

}

 Der Beweis verl\"auft analog zu dem von Satz 1.

 \vskip+0.5cm  
  
Der Autor dankt  Prof. O.N. Deutscher  f\"{u}r  die fruchtbare Diskussion
und Prof. C. Elsholtz f\"ur  sprachlische Hinweise.

\vskip+0.5cm
{\bf Literaturverzeichnis}
\vskip+0.5cm

\noindent[1]\,
B. J. Birch, H. Davenport,
\,\,
{\it
Quadratic equations in several variables},
Mathematical Proceedings of the Cambridge Philosophical Society, {\bf 54}:2 (1958), 135
- 138.

\noindent [2]\,
J. W. S. Cassels,\,\,{\it
Bounds for the least solutions of homogeneous quadratic equations},
Mathematical Proceedings of the Cambridge Philosophical Society, {\bf 51}:2 (1955), 262
- 264.

\noindent [3]\,
J. W. S. Cassels,\,\,{\it
Rational Quadratic Forms}, Academic Press, 1978.

  \noindent [4]
  H. Davenport,\,\,
  {\it
   Note on a theorem of Cassels},
   Mathematical Proceedings of the Cambridge Philosophical Society, {\bf 53}:2 (1957),  539
- 540.

 \noindent [5]
 L. Fukshansky,\,\,
{\it Heights and quadratic forms: Cassels' theorem and its generalizations},
in {\it
Diophantine methods, lttices,
and arithmetic theory of quadratic forms},
Contemp. Math., {\bf 587},  Ammer. Math. Soc., Providence, RI, 2013, 77 - 94.
}

\noindent
[6]\,
 L. Fishman,
D. Kleinbock, K. Merrill, D. Simmons,\,\,{\it
Intrinsic Diophantine approximation on manifolds},
preprint verf\"ugbar unter  arXiv:1405.7650v4[math.NT] 24 Sep 2015.

\noindent
[7]\,
D. Kleinbock, K. Merrill,\,\,{\it
Rational approximation on Spheres},
preprint verf\"ugbar unter arXiv:1301.0989v4 [math.NT] 25 May 2013,
erscheint in Israel Math. J.

\noindent [8]\,
N. Moshchevitin,\,\,{\it
\"Uber die rationalen Punkte auf der Sph\"are}, erscheint in 
Monatshefte f\"ur Mathematik (2015), DOI 10.1007/s00605-015-0818-4.

\noindent [9]\,
H.P. Schlickewei,\,\,
{\it
Kleine Nullstellen homogener quadratischer Gleichungen
},
Monatshefte f\"ur Mathematik, {\bf 100} (1985), 35 - 45.

\noindent [10]\,
H.P. Schlickewei, W.M. Schmidt,\,\,
{\it
Quadratic Forms Which Have Only Large Zeros I},
Monatshefte f\"ur Mathematik, {\bf 105} (1988), 295 - 311.

\noindent [11]\,
H.P. Schlickewei, W.M. Schmidt,\,\,
{\it
Isotrope Unterr\"aume rationaler quadratischer Formen},
Mathematische Zeitschrift, {\bf 201} (1989), 191 - 208.

\noindent [12]
R. Schulze-Pillot, {\it Small linearly independent zeros of quadratic forms},
Monatshefte f\"ur Mathematik,
{\bf  95} (1983),  24 -
249.

\end{document}